\documentclass[12pt,a4paper]{article}
\usepackage{amsmath,amsfonts,amssymb,amsthm,}
\usepackage[english]{babel}

\newtheorem{theorem}{Theorem}

\newtheorem{theor}{Theorem}

\theoremstyle{definition}

\theoremstyle{remark}

\title{Carleson measures on simply connected domains}

\author{ Mar\'ia\ J.\ Gonz\'alez}

\date{}

\begin{document}
	
	\maketitle
	
	\baselineskip=6mm
	\parskip=2.5mm
	\let\thefootnote\relax\footnote {\textit{2010 Mathematics Subject Classification:} 30C62, 42B20, 42B35}
	\let\thefootnote\relax\footnote{\textit{Keywords and phrases:}  Carleson measures, Hardy and Bergman spaces, BMO domains, chord-arc domains.}
	\begin{abstract}
		We study the Carleson measures associated to the Hardy and weighted Bergman spaces defined on general simply connected domains. This program was initiated by Zinsmeister in his paper \textit{Les domaines de Carleson }(1989), where he shows that the geometry of the domain plays a fundamental role. We will review the classical results presenting, in some cases, alternative proofs and will examine the situation for the weighted Bergman spaces.
	\end{abstract}

	\section*{Introduction}

	Let $\mathbb{D}$ denote the unit disc $\{z\in \mathbb{C};~|z|<1\}$ and let $\mathbb {T}=\partial \mathbb{D}$. If $I\subset \mathbb{T}$ is an interval, the Carleson square $S(I)\subset \mathbb{D}$ is the set
$$
S(I) = \{re^{i\theta} : e^{i\theta}\in I, 1-\frac{|I|}{2\pi}\leq r<1\}
$$
where $|I|$ denotes the length  of the interval $I$.

We will consider two classical spaces of analytic functions in $\mathbb{D}$: the Hardy spaces  and the weighted Bergman spaces. An analytic function $f$ in $\mathbb{D}$ is in the Hardy space $\mathcal{H}^p,~ 0< p <\infty $, if
$$
\sup_{0<r<1} ~\frac{1}{2\pi} \int_0^{2\pi}|f (re^{i\theta})|^p ~d\theta = \|f\|_{\mathcal{H}^p}^p < \infty
$$

 It is well known that  a function  $f\in \mathcal{H}^p$  has almost everywhere  non-tangential boundary limit  $f(e^{i\theta})\in L^p (\mathbb {T})$, and $ \|f\|_{\mathcal{H}^p}^p = \frac{1}{2\pi} \int_0^{2\pi}|f (e^{i\theta})|^p ~d\theta$.

The  weighted Bergman spaces $\mathcal{A}^p_\alpha$, $0< p <\infty, ~\alpha > -1$ consists of those analytic functions $f$ in $\mathbb{D}$ such that
$$
\frac{1}{\pi}\int_{\mathbb{D}}|f(z)|^p (1-|z|)^\alpha ~dm(z)= \|f\|_{\mathcal{A}^p_\alpha}^p < \infty
$$
 \noindent where $dm$ represents the Lebesgue area measure. When $\alpha=0$ we obtain the classical Bergman spaces.
  We refer to \cite{D}, \cite {G} and \cite {HKZ} for the theory on Hardy and Bergman spaces.

  A remarkable theorem due to   Carleson   \cite{C} states that a positive measure $\mu$ defined in $\mathbb{D}$ satisfies that for all $0<p<\infty$,
	$$
\int_{\mathbb{D}}|f(z)|^p~d\mu(z) \leq C_p ~\|f\|_{\mathcal{H}^p}^p,~ \textrm{for all}~ f\in \mathcal{H}^p
	$$
	if and only if there is a constant $c>0$ such that for all intervals $I\subset \mathbb{T} $, $\mu (S(I))< c |I|$.

	This theorem has led to various generalizations, including analogous results for the weighted Bergman spaces $\mathcal{A}^p_\alpha$  which can be summarized in the following theorem:
	
	\begin{theorem}[see \cite{OP,H,S}]\label{theoA}
		Let $0<p\leq q<\infty, ~\alpha>-1$ and let $\mu$ be a positive measure on $\mathbb{D}$. Then there exists  $C>0$ such that
		\begin{equation}\label{lue}
		\bigg(\int_{\mathbb{D}}|f(z)|^q  ~d\mu(z)\bigg)^{1/q}\leq C~ \|f\|_{\mathcal{A}^p_\alpha}
		\end{equation}
		\noindent if and only if  $\mu (S(I))\leq c~|I|^{(2+\alpha)q/p}$, for some $c>0$ and  for all intervals $I\subset \mathbb{T}$.
	\end{theorem}

 It has been pointed out by Luecking ( Th. 2.2 in \cite{Lu}) that a similar characterization can be expressed in terms of pseudohyperbolic balls. Indeed, it is easy to show that when the exponent $\beta >1$, the condition on the measure $\mu (S(I))< c~|I|^\beta$ for all intervals $I\subset \mathbb{T}$ is equivalent to the condition $\mu (B(z, r))< Cr^\beta$, for all $ z\in \mathbb{D}, ~r =\frac {1}{2} \operatorname{dist}(z,\partial \mathbb{D})$.

	Theorem A was obtained by Oleinik and Pavlov \cite{OP}, and independently by Hastings \cite{H} for $\alpha=0$  and Stegenga \cite{S} when $ p = q > 1$. Extensions of these results, where derivatives of the  functions are considered on left inside of (\ref{lue}) can be found in \cite{Lu} and the references within.
	
	In the context of Hardy spaces, Duren [D] extended Carleson's result to the range of exponents $q>p>0$ obtaining a similar  condition on the measure as the one given by Carleson in the case $q=p$.
	
	\begin{theorem}[see \cite{C,D}]\label{theoB}
		Let $0<p\leq q<\infty$ and let $\mu$ be a positive measure on $\mathbb{D}$. Then there exists  $C>0$ such that
		$$
		\bigg(\int_{\mathbb{D}}|f(z)|^q  ~d\mu(z)\bigg)^{1/q}\leq C~ \|f\|_{\mathcal{H}^p}
		$$
		
		\noindent if and only if  for some $c>0$,   $\mu (S(I))\leq c~|I|^{q/p}$, for all intervals $I\subset \mathbb{T}$.
		
	\end{theorem}

We will give a different proof of the sufficiency condition in Duren's result, which explicitly will show that it is enough to consider  pseudohyperbolic  balls instead of Carleson squares, although as we mentioned above they are equivalent since the exponent $q/p>1$.
	
	We have avoided up to this point  the term \textit{Carleson measures} because, depending on the authors,  the term is either used to define the measures for which the space $L^p(\mu)$ embeds continuously on the space of analytic functions under consideration or, in some other contexts, it is used to define the geometric characterization of the measure in terms of Carleson squares. The theorems above show that, for the Hardy and for the weighted Bergman spaces in the disc, both definitions are actually equivalent, but as we will see,  this will not be the case in a more general setting. We will adopt the following terminology:

Given a  Banach space $X$ of analytic functions in a domain $\Omega\subset \mathbb{C}$ with norm $ \|.\|$, we say that a positive measure $\mu$ in $\Omega$
is a \textit{$q$-Carleson measure for X}  if there exists a
constant $C>0$ such that
$$
\bigg(\int_{\mathbb{D}}|f(z)|^q  ~d\mu(z)\bigg)^{1/q}\leq C~ \|f\|_X
$$

	We are interested in characterizing the Carleson measures for spaces of analytic functions defined on bounded simply connected domains $\Omega \subset \mathbb{C}$. This problem has been initially studied by Zinsmeister in \cite{Z} where he extended Carleson`s result to more general domains, showing how the geometry of the domain plays a fundamental role.

	 In order to define the Hardy spaces $\mathcal{H}^p (\Omega)$ we need to assume initially that $\partial \Omega$ is  locally rectifiable. Denoting by $ds$ the arc length measure, we define for $0<p<\infty$
	$$
	\mathcal{H}^p (\Omega)= \{f ~\text{analytic in}~ \Omega ~; ~ \int_{\partial\Omega}|f(z)|^p~ds < \infty\}
	$$
	
The theory of  Hardy spaces is well understood when $\Omega$ is a chord-arc domain, that is a domain bounded by a chord arc curve. In this case, the functions in  $\mathcal{H}^p (\Omega)$  can be characterized in terms of the nontangential maximal function and the area integral as in the classical case, (see \cite{JK}). Recall that a   locally rectifiable curve $\Gamma$  is a chord-arc curve if $\ell_{\Gamma}(z_1,z_2)\le
		K|z_1-z_2|$ for some $K>0$ and for all $z_1,z_2\in\Gamma$, where $\ell_{\Gamma}(z_1,z_2)$ denotes the length of the  shortest arc
		of $\Gamma$ joining $z_1$ and $z_2$.

Let $\Omega$ be a simply connected domain and  $\varphi: \mathbb{D}\rightarrow \Omega$  a conformal map. Let $\mu$ be  a positive measure on $\Omega$. When $\partial \Omega$ is  rectifiable, a simple change of variables  shows that $\mu$ being a $q$-Carleson measure for $\mathcal{H}^p (\Omega)$, that is,		
$$
		\bigg(\int_\Omega|f|^q  ~d\mu\bigg)^{1/q}\leq C~ \bigg(\int_{\partial\Omega}|f|^p  ~ds\bigg)^{1/p}
		$$
is equivalent to
$$
		\bigg(\int_{\mathbb{D}}|f\circ \varphi|^q  ~d\varphi^*(\mu)\bigg)^{1/q}\leq C~ \bigg(\int_{\partial\mathbb{D}}|f|^p |\varphi'| ~ds\bigg)^{1/p}
		$$
\noindent where  $\varphi^*(\mu)$ denotes the pullback of~$\mu$, that is $\varphi^*(\mu) (E)=\mu (\varphi(E))$, for any set $ E\subset \mathbb{D}$.

Therefore, as in \cite {Z}, if we define  the measure $\nu$ in $ \mathbb{D}$ as $\nu=\frac{1}{|\varphi'|^{q/p}} ~\varphi^*(\mu)$,  we obtain the following observation that we state as a remark for further references.

\noindent \textbf{Remark:}
 The measure $\mu$ is a  $q$-Carleson measure for $\mathcal{H}^p(\Omega)$ if and only if the measure $\nu$ is  $q$-Carleson measure for $\mathcal{H}^p (\mathbb{D})$

 The advantage of this point of view is that we do not need to assume rectifiability on the boundary of the domain in order to study  Carleson measures for $\mathcal{H}^p(\Omega)$.

	\begin{theorem}[\cite{Z}]\label{theoC}
		Let $\Omega$ be a bounded simply connected domain and let  $\varphi: \mathbb{D}\rightarrow \Omega$ be a conformal map. Assume that $\log{ \varphi '}\in \text{BMOA}(\mathbb{D})$. Then a positive measure $\mu$ in $\Omega$ is a $p$-Carleson measure  for  $\mathcal{H}^p (\Omega), 0<p<\infty$, if for some $c>0$,
		$ \mu (B(\xi, R)\cap \Omega)\leq c R$, for all $\xi\in \partial \Omega $ and $R>0$.
		\end {theorem}
Moreover, it is also showed  in \cite {Z} that the condition on the domain: $\log{ \varphi '}\in \text{BMOA}(\mathbb{D})$, is a necessary condition for Theorem C to hold.

The geometry of domains for which $\log{ \varphi '}\in \text{BMOA}(\mathbb{D})$ has been characterized by Bishop and Jones in \cite {BJ}. The boundary of these domains might not be rectifiable, and though we will not give the precise definition, let us just mention that in  some sense they are rectifiable most of the time on all scales. A typical example is a variant of the snowflake where at each iteration step, one of sides
of the triangle, for instance the left one, is left unchanged.

To prove theorem C in \cite{Z}, it is first shown  that the result holds for chord arc domains. The general result follows using the fact that when the BMO norm is small enough the domain is chord arc. We will give a different proof based on a stopping time argument.

		For the converse result of theorem C one needs a stronger assumption on the boundary of the domain. We say that a curve $\Gamma$ is  Ahlfors regular if there exists  $C>0$ such that for all
$z_0 \in\Gamma$ and all $R > 0$, the arclength of $B(z_0, R) \cap\Gamma\leq C R$. These curves were studied by G. David in the context of the Cauchy integral \cite{Da}.

Chord arc curves are Ahlfors regular but not viceversa, for example cusps are Ahlfors regular but not chord arc. On the other hand, if we add the condition that the curve is a quasicircle, then  Ahlfors regular and chord arc are equivalent.

		\begin{theorem}[\cite{Z}]\label{theoD}
			Let $\Omega$ be a bounded simply connected domain, and $\mu$ be a positive measure in $\Omega$. Assume that $\partial\Omega$ is an Ahlfors-regular Jordan curve. Then a positive measure $\mu$ in $\Omega$ is a $p$-Carleson measure  for  $\mathcal{H}^p (\Omega), 0<p<\infty$,  if and only if
	for some $c>0$,	$ \mu (B(\xi, R)\cap \Omega)\leq c R$, for all $\xi\in \partial \Omega $ and $R>0$.
		\end{theorem}

The sufficiency condition follows immediately from Theorem C, since domains bounded by Ahlfors regular curves satisfy that $\log{ \varphi '}\in \text{BMOA}(\mathbb{D})$. The necessity is proved using Heyman-Wu Theorem. We will provide a simple proof in the case that $\partial\Omega$ is also a quasicircle, and therefore chord-arc.

We will state the  next results in terms of Whitney balls in the domain $\Omega$. They play the same role as  pseudohyperbolic  balls in the unit disc. We will say that a ball $B(z,r)\subset \Omega$ is a Whitney ball if $c ~\delta_\Omega (z)\leq r\leq 1/2 ~\delta_\Omega (z)$, for some fixed constant $c>0$, where $\delta_\Omega (z)$ denotes the distance from $z$ to the boundary of $\Omega$.

In the context of Hardy spaces  we obtain the analogous result of Duren's theorem in the classical case.

\begin{theor}
		Let $\Omega$ be a bounded simply connected domain,and $\mu$ a positive measure in $\Omega$. If $0<p<q<\infty$, then $\mu$ is a $q$-Carleson measure
		for  $\mathcal{H}^p (\Omega)$ if and only if there is a constant $c>0$ such that  $ \mu (B(z, r))\leq c r^{q/p}$, for all Whitney balls $B(z,r)\subset \Omega$.
		\end {theor}

Notice that in contrast with the classical case $p=q$, no condition on the geometry of the domain is assumed. This is not surprising, it is a consequence that analogously as what happens in the disc, when the exponents are bigger that 1, the characterizations of the Carleson measures  can be expressed in terms of Whitney balls instead of balls centered at the boundary.

 We define the Bergman spaces $\mathcal{A}^p_\alpha (\Omega)$ for $0<p<\infty$, $\alpha>-1$ as
	$$
	\mathcal{A}^p_\alpha (\Omega)= \{f ~\text{analytic in}~ \Omega ~; ~ \int_{\Omega}|f(z)|^p ~ ~\delta_{ \Omega} (z)^\alpha ~dm(z)< \infty\}
	$$

In this setting, we will  not only consider analytic functions but  quasi-subharmonic functions as well. Let $u \geq 0$ be a locally bounded, measurable function on $\Omega$. We say
	that the function $u$ is \textit{$C$-quasi-nearly subharmonic}  if the following condition
	is satisfied:
	\begin{equation}\label{subh}
	u(a)\leq C \frac{1}{r^2}\int_{B(a,r)} u~dm
	\end{equation}
	\noindent  whenever $B(a, r)\subset \Omega$.
	
	\noindent One can view (\ref{subh}) as a weak mean value property. Besides of nonnegative
	subharmonic functions, it also holds for nonnegative powers of subharmonic functions, and for subsolutions to a large family of second order elliptic equations, see (\cite {P}, \cite{PR}).

	\begin{theor} Let $\Omega$ be a bounded simply connected domain  and $\mu$ be a positive measure on $\Omega$. Let $0<p\leq q<\infty$, and $\alpha>-1$.
		 The following conditions  are equivalent:
		\begin{enumerate}
			\item [(i)]
			There is $c>0$ such that $\mu (B(z, r))\leq c ~r^{(2+\alpha)q/p}$, for all Whitney balls $B(z,r)\subset \Omega$.
			\item [(ii)]
			There is $C'>0$ such that for any $C$-quasi-subharmonic function $g\geq 0$ in $\Omega$
			$$
			\bigg(\int_\Omega g^q(z)  ~d\mu(z)\bigg)^{1/q}\leq C'~ \bigg(\int_{\Omega}g^p(z)  ~\delta_{ \Omega} (z)^\alpha ~dm(z)\bigg)^{1/p}
			$$
			
			\item [(iii)]
			The measure $\mu$ is a $q$-Carleson measure for the space $\mathcal{A}^p_\alpha (\Omega)$.
			
		\end{enumerate}
	\end{theor}
	
The proof of (i) implies (ii) is very similar to the one given by Luecking in \cite{L} where he proves that the analogous statement holds for subharmonic functions in the unit disc.

The structure of the paper is as follows: In section 1 we will fix the notation, definitions and state some basic results. The proofs of the already known results will be shown in section 2, that is, the sufficiency condition in Duren's theorem, Theorem C, and the chord-arc case in Theorem D. Finally, in Section 3, we will prove Theorem 1 and Theorem 2.

		\section{Basic facts and definitions}\label{sec1}
		
In the paper, the letter $C$ denotes a constant that may change at different
		occurrences. The notation $A \simeq B$ means that there is a constant $C$  such that $1/C.A\leq B \leq C.A$. The notation $A\lesssim  B$ ( $A\gtrsim B$)
		means that there is a constant $C$  such that $A \leq C.B$ ($A\geq C.B$).

 Also, as usual, we denote by $\mathbb{T}$  the boundary of the unit disc, and $B(z_0,R)$ the ball of radius $R$ centered at the point $z_0\in\mathbb{C}$. If $B$ is a ball, $2B$ is the ball with
		the same center as $B$ and twice the radius of $B$, and similarly for squares. Given a domain $\Omega\in \mathbb{C}$, for any  $z\in \Omega$  we set $\delta_\Omega (z)=\operatorname{dist} (z,\partial\Omega)$. If the context is clear we will drop the subindex $\Omega$ and simply write $\delta (z)$.

Given an interval $I\in \mathbb{T}$, and the corresponding Carleson square $S(I)$, we define the top of the square $T(S)$ as the set of points
$$
T(S) = \{re^{i\theta} : e^{i\theta}\in I, 1-\frac{|I|}{2\pi}\leq r<  1-\frac{|I|}{4\pi}\}
$$
		
		 A locally integrable function~$f$ belongs to the
		space~$\operatorname{BMO}(\mathbb{T})$ if
		$$
		\|f\|_*=\sup_I\frac{1}{|I|} \int_I |f(x)-a_I|\,dx<\infty
		$$
		where the supremum is taken over all arcs $I\in \mathbb{T}$ and $a_I=\frac{1}{|I|}\int_I f(y)\,dy$

It is a well known result, see for instance Th.1.2, Ch.VI in \cite{G}, that if $f\in \operatorname{BMO}(\mathbb{T})$, then
		\begin{equation*}
        \sup_{z\in \mathbb{D}}\int_ {\mathbb{T}}|f(\xi)-f(z)| ~P_z(\xi) ~|d\xi| =A <\infty
        \end{equation*}
where $f(z)= \int_{\mathbb{T}} P_z (\xi) f(\xi) ~|d\xi| $ is the Poisson integral of $f$. Moreover $A \simeq |f\|_*$.

		  In particular, for any interval $I\in \mathbb{T}$, if $z_I$ denotes the point $z_I= (1-|I|/2)\xi_I$, being $\xi_I$ the midpoint of the interval $I$, we have
\begin{equation}\label{BMOA}
     \frac{1}{|I|} \int_I |f(\xi)-f(z)| ~d\xi \leq c \|f\|_*
     \end{equation}
for some $c=c(\|f\|_*)$.

The space $\operatorname{BMOA}(\mathbb{D})$ consists of those functions in the Hardy space $\mathcal{H}^1(\mathbb{D})$ whose boundary values are in $\operatorname{BMO}(\mathbb{T})$. We refer to \cite{G} and \cite{Po} for the main properties of $\operatorname{BMOA}$ functions.

For a harmonic function $f$ in $\mathbb{D}$ and for any $~\xi\in\partial \mathbb{D}$, the nontangential maximal function
$f^\ast$ is defined as $f^\ast(\xi)= \operatorname{sup} \{|f(z)|: z\in\Gamma_\xi\}$, where $\Gamma_\xi$ denotes the cone
$\Gamma_\xi = \{z \in \mathbb{D} : |z - \xi| < \alpha (1-|z|)\}$, for some fixed $\alpha>0$.

Given a function $f\in L^1_{\operatorname{loc}} (\mathbb{R})$, the Hardy-Littlewood maximal function of $f$ is
$$
Mf (x) = \sup_{x\in I}\frac{1}{|I|}\int_I |f(t)| ~dt
$$

It is well known that the operator $Mf$ is bounded in $L^p; 1< p< \infty$ and it is of weak type $1-1$, that is
$$
|{t \in \mathbb{R}: Mf (t) >\lambda }| \leq \frac{2}{\lambda}\|f\|_1, ~ \lambda >0.
$$

The importance of the maximal function $Mf$ is that it
majorizes many other functions associated with $f$, such as  the non-tangential maximal function of the Poisson integral of f. Therefore, if $u(z)= P_z\ast f$, its  non-tangential maximal function $u^\ast(x)=\sup_{\Gamma_x}|u(z)|$ satisfies
\begin{equation}\label{maximal}
|x\in I; ~u^\ast(x)> \lambda|\leq  \frac{c}{\lambda}\int_I |f(x)|~dx
\end{equation}
\noindent for some constant $c>0$ depending only on the aperture of the cone $\Gamma_x$, (see for instance Chapter I.4 in \cite{G} for these and many others related results).

We finish this section by mentioning some well known results on conformal mappings, we refer the reader to Ch.I in \cite{Po} for an excellent overview.

Simply connected, proper subdomains of the plane inherit a hyperbolic metric
from the unit disk via the Riemann map. If $\varphi:\Omega\rightarrow \mathbb{D}$ is conformal and $w =\varphi(z)$
then $\rho_\Omega(w_1, w_2) = \rho_\mathbb{D}(z_1, z_2)$ defines the hyperbolic metric on $\Omega$ and is independent
of the particular choice of $\varphi$. It is often convenient to estimate $\rho_\Omega$ in terms of the
more geometric quasi-hyperbolic metric on $\Omega$ which is defined as
$$
\tilde{\rho}_\Omega(w_1, w_2) = \operatorname{inf}\int_{w_1}^{w_2}\frac{|dw|}{\delta_\Omega(w)}
$$
where the infimum is taken over all arcs in $\Omega$ joining $w_1$ to $w_2$. It follows from Koebe $1/4$
theorem that the two metrics are comparable with bounds independent of the domain.
A Whitney decomposition of the domain $\Omega$ is a covering of $\Omega$ by squares $Q_k$ with
disjoint interiors and the property that $\operatorname {diam}(Q_k) \simeq \delta_\Omega (Q_k)$. By our remarks
above, each square in a Whitney decomposition has uniformly bounded hyperbolic
diameter (and contains a ball with hyperbolic radius bounded uniformly from below). Thus bounding the hyperbolic length of a path often reduces to simply estimating
the number of Whitney squares it hits.

 Analogously, we will say that a ball $B(z,r)\in \Omega$ is a Whitney ball if $c ~\delta_\Omega (z)\leq r\leq 1/2 ~\delta_\Omega (z)$, for some fixed constant $c>0$.

 Let us also recall that the function $\log{\varphi}'\in \mathcal{B}$, where $\mathcal{B}$ denotes the Bloch space in $\mathbb{D}$. Therefore, for any Carleson square $S\in \mathbb{D}$, if $z_1,z_2$ are any two points in the top of the square $ T(S)$,  then
 \begin{equation}\label{conformal}
 |\varphi'(z_1)|\simeq |\varphi'(z_2)|;~z_1,z_2\in T(S)
 \end{equation}
  This is because the hyperbolic diameter of the top of the squares is uniformly bounded.

		\section{Proofs of some known results}\label{sec2}

We begin this section by giving a short proof of the characterizations of the $q$-Carleson measures for $\mathcal{H}^p; ~q>p>0$ due to Duren.		
		
\begin{proof}{\textbf{(Proof of Theorem B, $q>p>0$)}}
\quad

 It is enough to prove the result for  that $p=1$, so let $q>1$. Assume first that for any Whitney ball $B(z, r)\subset \mathbb{D}$
\begin{equation}\label{medida}
\mu (B(z, r))\leq c~r^q
\end{equation}

\noindent It is a well known result  that if $f\in \mathcal{H}^1$, then $f(z)\lesssim 1/(1-|z|)$ for $ z\in \mathbb{D}$,  and the non-tangential maximal function $f^\ast\in \mathcal{H}^1$, see for example \cite{G}.
 Using these results and Fubini's theorem, we get
 \begin{equation*}
 \begin{split}
 \int_\mathbb{D}|f(z)|^q~d\mu(z)=\int_{\partial \mathbb{D}} \int_{\Gamma_\xi}\frac{|f(z)|^q}{1-|z|}~d\mu(z)~|d\xi|\leq
 \int_{\partial \mathbb{D}}f^\ast(\xi) \int_{\Gamma_\xi}\frac{|f(z)|^{q-1}}{1-|z|}~d\mu(z)~|d\xi|\\
 \lesssim \int_{\partial \mathbb{D}}f^\ast(\xi) \int_{\Gamma_\xi}\frac{1}{(1-|z|)^{2-q}}~d\mu(z)~|d\xi|\lesssim \int_{\partial \mathbb{D}}f^\ast(\xi) ~|d\xi|\lesssim ||f||_{\mathcal{H}^1}
\end{split}
\end{equation*}
\noindent since the integral on the cone $\Gamma_\xi$ can be estimated in terms of Whitney balls $B_k$ centered at points $z_k=(1-2^{-k})\xi;~k=0,1,... $ as
$$
 \int_{\Gamma_\xi}\frac{1}{(1-|z|)^{2-q}}~d\mu(z)\lesssim\sum_k \mu(B_k)/(1-|z_k|)^{2-q}\lesssim  \sum_k (2^{-k})^{2q-2}\lesssim 1
$$
\noindent because $\mu$ satisfies (\ref{medida}) and $q>1$

The converse result is standard, we give the idea of the proof for the sake of completeness. Assume that $\mu$ is a $q$-Carleson measure for $\mathcal{H}^1$.
 For each $z_0\in\mathbb{D}$, choose
 $z_0^\ast \in \mathbb{C}\setminus \mathbb{D}$ such that $|z_0-z_0^\ast| \sim 1-|z_0| \sim 1-|z_0^\ast|$. Then it is easy to prove that $f(z)=\frac{1}{(z-z_0^\ast)^2}$ is in $\mathcal{H}^1$ with $\|f\|_{\mathcal{H}^1}\sim 1/(1-|z_0|)$. The result now easily follows by applying the hypothesis to the function $f$.
\end{proof}

Next we prove Zinsmeister's result in Theorem C. Thus assuming that $\log{\varphi'}\in
\operatorname{BMOA}$ we want to show that, if $\mu(B(\xi,r)\cap\Omega)\leq Cr; ~\xi\in \partial \Omega$ then $\mu$ is a Carleson measure for the Hardy spaces in $\Omega$ or, equivalently by the remark in the Introduction, that  $\nu(B(\xi,r)\cap\mathbb{D} )\leq c~r; ~\xi\in \partial \mathbb{D},~r>0$.

To simplify the notation we will replace the ball $B(\xi,r); ~\xi\in \partial \Omega$  by a Carleson square in the upper half plane, that is, a square with base on some interval $ I\subset \mathbb{R}$. So, let $ Q= \{(x,y); x\in I, 0 \leq y\leq |I|\}$ . For any such square, we define the top of the square $T(Q)=\{(x,y); x\in I, 1/2 \leq y\leq |I|\}$, and the center of $T(Q)$ as the point $z_I=x_I+i3/2 |I|$, where $x_I$ denotes the midpoint of the interval $I$.

\begin{proof}{\textbf{(Proof of Theorem C)}}
\quad

Let $I$ be any interval in $ [0,1]$, define the functions
$f_I(z)= \log{\varphi'(z)}- \log{\varphi'(z_I)}$
and  $u_I(z)=\operatorname {Re}f(z)= \log{|\varphi'(z)|}- \log{|\varphi'(z_I)|}$. Since $f_I\in\mathcal{H}^1$, the harmonic function $u_I$ is the Poisson integral of its boundary values, i.e. $u_I(z)= P_z\ast (\log{|\varphi'(x)|}- \log{|\varphi'(z_I)}|)$ (see Th.3.6, Ch.II in\cite{G}).
On the other hand, by (\ref{BMOA}) we get
$$
\frac {1}{|I|} \int_I u_I(x) ~dx =\frac {1}{|I|} \int_I\big|\log{|\varphi'(x)|}- \log{|\varphi'(z_I)|}\big|~dx
\lesssim  \|\log{\varphi'}\|_\ast
$$
Therefore by (\ref{maximal}), for any interval $I$
\begin{equation}\label{Tch}
\frac{1}{|I|}|x\in I; ~u_I^\ast (x)> \lambda|\lesssim \frac{ \|\log{\varphi'(z)}\|_\ast}{\lambda}
\end{equation}

Fix now an interval $I$ and let $Q$ be the  corresponding Carleson square. We want to show that $\nu(Q)\leq c |I|$.

 The idea  is to divide  $Q$ into a countable union of disjoint regions, in such a way that $|\varphi'|$ behaves like a constant inside  each of those regions. For that, we will use a stopping time argument.

 For each $k=1,2,...$ form the $2^k$ intervals obtained by dividing $I$ dyadically, and associate to each of these intervals the corresponding Carleson square. We obtained in this way $2^k$ squares of length $2^{-k} |I|$. Denote by $\{Q_j\}$ these collection of squares. Note that each $Q_j$ is associated to some $k\in \mathbb{N}$ and it is contained in some other  square $Q_l$ associated to $k-1$, which we will call the father ( if $k=1$, the father is $Q$). Squares with the same father will be called brothers.

 Let $M>0$ be a big enough constant that will be fixed later. Recall that $z_I$ is the center of $T(Q)$. Define the first generation $G_1(Q)$ as those $Q_j \subset Q$ such that
 \begin{equation}\label{stop}
 \sup_{z\in T(Q_j)} |\log|{\varphi'(z)|}- \log |{\varphi'(z_I)|}| > \log{M}
 \end{equation}
 \noindent and  $Q_j $ is maximal.

  We also say that the interval $I_j\in G_1(I)$, if $I_j$ is the base of some square  $Q_j\in G_1(Q)$. Note that the intervals $\{I_j\}; ~I_j\in G_1(I)$ have pairwise disjoint interiors, and that in the region $\mathcal{R}_1(Q)=Q\setminus \bigcup_{G_1(Q)}Q_j$,
 \begin{equation}\label{bilip}
 \frac{1}{M}|\varphi'(z_I)|\leq  |\varphi'(z)|\leq M |\varphi'(z_I);~z\in \mathcal{R}_1
 \end{equation}

  On the other hand, if $x\in I_j$ for some $I_j\in G_1(I)$, its cone $\Gamma_x$ contains points which belong to $T(Q_j)$. So, by the choice of the squares in the first generation and by (\ref{conformal}), we get that $u_I^\ast(x)> \log{c M}$, for some universal constant $c>0$. Therefore, by (\ref{Tch}),
 \begin{equation}\label{intervals}
 \sum_{I_j\in G_1(I)}|I_j|\leq |x\in I; ~u_I^\ast (x)> c M|\leq \frac{1}{10} |I|
 \end{equation}
\noindent  if $M> M_0$ where $M_0$ is a constant depending only on $\|\log{\varphi'(z)}\|_\ast$.

Next, to each square in the first generation $Q_j\in G_1(Q)$ we apply the same stopping time rule (\ref{stop}), that is we start partitioning $Q_j$ and we choose those maximal squares $Q_k\subset Q_j$ for which
 \begin{equation*}
 \sup_{z\in T(Q_k)} |\log|{\varphi'(z)|}- \log |{\varphi'(z_{I_j})|}| > \log{M}
 \end{equation*}
 where $z_{I_j}$ is now the center of $T(Q_j)$.  In this way we obtain  a new generation of squares $G_2(Q_j)$. Define the second generation as
$$
G_2(Q)=\cup\{G_1(Q_j^1);~Q_j^1\subset G_1\}=\{Q_1^2,Q_2^2,...\}
$$

Repeating  the process with $G_2$ and continuing inductively we obtain later generations of squares $G_n(Q)= \{Q_1^n,Q_2^n,...\}$, and generations of intervals $G_n(I)= \{I_1^n,I_2^n,...\},~n\in \mathcal{N}$. Moreover, by (\ref{intervals}),
\begin{equation}\label{int}
\sum_{I_j\in G_n(I)}|I_j|\leq  \frac{1}{10} \sum_{I_j\in G_{n-1}(I)}|I_j|\leq....\leq  \frac{1}{10^n} |I|
\end{equation}

We are ready now to prove the theorem. Let $Q$ be the initial Carleson square. By the definition of the measure $\nu$
\begin{equation}\label{est}
\nu(Q)=\int_Q \frac{1}{|\varphi'(z)|}~d\mu(z)= \sum_{n=1} \sum_{l;\mathcal{R}_l\in G_n}\int_{\mathcal{R}_l}\frac{1}{|\varphi'(z)|}~d\mu(z)
\end{equation}

 By construction, the region $\mathcal{R}_l\in G_n$  is contained in some square  of the previous generation, that we will denote by $Q_l$. As it was observed in (\ref{bilip}), if $z\in \mathcal{R}_l$ then $|\varphi'(z)|\simeq |\varphi'(z_l)|$, where $z_l$ is the center of $T(Q_l)$.

\noindent Therefore $|\varphi (z)-\varphi(w)|\simeq |\varphi'(z_l)||z-w|$ for all$ ~ z,w \in \mathcal{R}_l$. Since the diameter of $R_l$ is comparable to $\operatorname{Im}z_l$,  Koebe's theorem
 implies that  $\varphi (\mathcal{R}_l)\subset B(\varphi (z_l),\delta _\Omega(z_l))\cap \Omega$. By the hypothesis on the measure $\mu$
$$
\mu (\varphi (\mathcal{R}_l))\leq \mu  (B(\varphi (z_l),\delta _\Omega(z_l))\cap \Omega)\lesssim \delta _\Omega(z_l)\simeq |\varphi'(z_l)|\operatorname{Im}z_l
$$
We conclude by (\ref{est}) and (\ref{int}) that
$$
\nu(Q)\simeq \sum_{n=1} \sum_{l;\mathcal{R}_l\in G_n}1/|\varphi'(z_l)|~\mu (\varphi (\mathcal{R}_l))\lesssim  \sum_{n=1} \sum_{l;\mathcal{R}_l\in G_n}\operatorname{Im}(z_l)\lesssim  \sum_{n=1} \sum_j |I_j^{n-1}|\lesssim \sum_{n=1} \frac{1}{10^n}|I|\leq c |I|
$$
\noindent with comparison constants $c=c(||\mu\|, \|\log{\varphi'}\|_\ast)$.
\end{proof}

We finish the section by giving an alternative proof  of the theorem by Zinsmeister, on the characterization of Carleson measures on domains bounded by Ahlfors regular curves, when  we add the extra hypothesis  that the curve is  a quasicircle, and therefore chord arc.

\begin{proof} (\textbf{Proof of Theorem D for chord-arc domains})

Let $\Omega$ be a domain bounded by a chord arc curve $\Gamma$.  Assume that $\mu$ satisfies  a $1$-Carleson measure for $\mathcal{H}^1(\Omega)$, that is, for all $f\in\mathcal{H}^1(\Omega)$,
\begin{equation} \label{reg}
\int_\Omega |f|~d\mu \lesssim \int_{\Gamma} |f| ~ds
\end{equation}

 Let $B(\xi_0,R)$ be a ball centered at $\xi_0\in \partial \Omega$ of radius $R>0$, and $w_0$ be a point in $B(\xi_0,R)\cap\Omega$ such that $\delta (w_0) \simeq R$. Since chord arc curves are quasicircles,  by the circular distortion theorem they admit a quasiconformal reflection (see \cite{A}). Thus   we can choose a point
 $w_0^\ast \in \mathbb{C}\setminus \Omega$ such that for all $w\in B(\xi_0,R)\cap\Omega$, it holds that  $|w-w_0^\ast| \sim \delta (w_0) \sim \delta (w_0^\ast)$.

  Consider the function  $f(z)=\frac{1}{(w-w_0^\ast)^2}$. It is easy to show that $f\in \mathcal{H}^1(\Omega)$ with $\|f\|_{\mathcal{H}^1}\lesssim 1/\delta (w_0)$. Indeed, consider the balls $B_k= B(w_0^\ast, 2^k \delta(w_0^\ast)); ~k=1,2,...$, and the annuli $A_k= B_k\setminus B_{k-1};~k=2,3,..$ Then
  $$
  \int_{\Gamma}\frac{1}{|w-w_0^\ast|^2}  ~ds = \sum_{k=2} \int_{ \Gamma \cap A_k}\frac{1}{|w-w_0^\ast|^2}~ds \lesssim \sum_{k=2} \frac{1}{(2^k \delta(w_0^\ast))^2} \operatorname{length} (\Gamma \cap B_k) \lesssim \frac{1}{\delta (w_0)}
 $$
  because $\Gamma$ is Ahlfors regular and therefore $ \operatorname{length} (\Gamma \cap B_k) \lesssim r(B_k)= 2^k \delta(w_0^\ast)$.

 We can then use (\ref{reg}) to bound $\mu(B(\xi_0,R)\cap \Omega)$ as follows
 $$
 \frac{\mu(B(\xi_0,R)\cap \Omega)}{\delta (w_0)^2}
 \simeq  \int_{B(\xi_0,R)\cap\Omega} \frac{1}{|w-w_0^\ast|^2}~d\mu(z)
  \lesssim\int_{\Gamma}\frac{1}{|w-w_0^\ast|^2}~ds \lesssim \frac{1}{\delta (w_0)}
 $$
 Hence, $\mu(B(\xi_0,R)\cap\Omega)\lesssim \delta(w_0)\simeq R$ as we wanted to show.	
 \end{proof}

\section{Proofs of  Theorems 1 and 2}\label{sec3}

\begin{proof}{\textbf{(Proof of Theorem 1)}}
			
\quad
Let $\varphi$ represent a conformal mapping from $\mathbb{D}$ onto $\Omega$. Assume first that the measure $\mu$ satisfies that
\begin{equation}\label {muu}
 \mu (B(w, r))\leq c ~r^{q/p}
 \end{equation}
 \noindent for all Whitney balls  $B(w, r)\subset\Omega$.

 Define the measure $\nu$ in $\mathbb{D}$ as $d\nu=\frac{1}{|\varphi'|^{q/p}}~d\varphi^\ast (\mu)$. By the remark in the Introduction, showing that $\mu$ is a Carleson measure in $\Omega$ is equivalent to showing that $\nu$ is a Carleson measure in $\mathbb{D}$.  So, let $B(z_0,r)\subset \mathbb{D}$ be a Whitney ball in $\mathbb{D}$. As $|\varphi'(z)|\simeq |\varphi'(z_0)|$ when $z\in B(z_0,r)$, we get
 $$
\nu(B(z_0,r))=\int _{B(z_0,r)} \frac{1}{|\varphi'(z)|^{q/p}}~d\varphi^\ast (\mu(z))\simeq \frac{1}{|\varphi'(z_0)|^{q/p}}\int_{\varphi(B(z_0,r))}~d\mu
$$

Set $w_0=\varphi(z_0)$. As explained in  section 1, $\varphi(B(z_0,r))$ can be covered by  a finite union of Whitney balls  $B(w_j, r_j)\subset \Omega;~j=1...,N_0$ where $N_0$ is a universal constant and  $r_j \simeq \delta (w_j)\simeq \delta (w_0)$. Therefore, by (\ref{muu}) and  by Koebe's theorem
\begin{equation*}
\begin{split}
\nu(B(z_0,r))\simeq \frac{1}{|\varphi'(z_0)|^{q/p}}~ \mu(\varphi(B(w_0,R)))\lesssim \frac{1}{|\varphi'(z_0)|^{q/p}}\sum_j \mu(B(w_j, r_j))\\
\leq \frac{1}{|\varphi'(z_0)|^{q/p}}\sum_j r_j^{q/p}\simeq \frac{1}{|\varphi'(z_0)|^{q/p}}\delta (w_0)^{q/p}\simeq (1-|z_0)^{q/p}\simeq r^{q/p}
\end{split}
\end{equation*}
\noindent which implies by Theorem B that $\nu$ is a Carleson measure in $\mathbb{D}$.

To prove the  "only if" part we proceed in a similar way. Assume now that $\mu$ is a $q$-Carleson measure for $\mathcal{H}^p (\Omega)$. This is equivalent, by the remark in the Introduction and  Duren's result in Theorem B, to the assumption that  $\nu (B(z,r))\lesssim r$ for all Whitney balls $B(z,r)\subset \mathbb{D}$.

 Next, consider a Whitney ball $B(w_0, R)\subset \Omega$ , i.e. $R\sim \delta(w_0)$, and let $\varphi^{-1}(w_0)=z_0 $. The set  $\varphi^{-1}(B(w_0,R))$  is contained in a finite union of Whitney balls  $B(z_j, r_j)\subset \mathbb{D};~j=1...,J_0$ for some universal constant $J_0$. Thus,  $|\varphi'(z)|\simeq |\varphi'(z_0)|$ when $z\in \cup_j B(z_j,r_j)$, and we obtain that
\begin{equation*}
\begin{split}
\mu(B(w_0,R)= \int_{B(w_0,R)}~d\mu = \int_{\varphi^{-1} (B(w_0,R))}  |\varphi'|^{q/p} ~d\nu \leq
\sum_j \int_{B(z_j,r_j)}|\varphi'|^{q/p} ~d\nu\\
\lesssim \sum_j |\varphi'(z_j)|^{q/p} ~ \nu(B(z_j,r_j))\lesssim \sum_j |\varphi'(z_0)|^{q/p}~ (1-|z_j|)^{q/p} \simeq  R^{q/p}
\end{split}
\end{equation*}

\noindent since $1-|z_j|\simeq 1-|z_0|$, and   by Koebe's theorem $\delta(w_0)\simeq \varphi'(z_0)~(1-|z_0|)$.

\end{proof}







\begin{proof}{\textbf{(Proof of Theorem 2)}}
			\quad

			$(i)\Rightarrow (ii)$ To simplify notation set $\delta (z)=\delta_\Omega (z)$. Let us begin with a simple observation related to Whitney balls: If $w\in B(z,\frac{1}{4}\delta (z))$, then $z\in B(w,\frac{1}{3}\delta (w))$. Besides $\delta (w)\simeq \delta (z)$.
			
			\noindent To see this, just note that $\delta (z) \leq |z-w|+\delta(w)$. Therefore, if $|w-z| \leq \frac{1}{4}~\delta (z)$, then $|w-z| \leq \frac{1}{4}(|z-w|+\delta (w))$, that is, $|w-z| \leq \frac{1}{3}~\delta (w)$. The rest of the statement follows by a similar argument.
			
			Let $g$ be a $C$-quasi-subharmonic function, then for any $p>0$, $g^p$ is $K$-quasi-subharmonic with constant $K=K(C,p)$ \cite{P}. Therefore, for all $z\in \Omega$ and any Whitney ball $B(z,r)$
			$$
			g^p(z)\lesssim \frac{1}{r^2}\int_{B(z,r)} g^p(w)~dm(w)
			$$
			We can then write

			\begin{equation}\label{ad}
		\Big(\int_\Omega g^q (z) ~d\mu(z) \Big)^{p/q}	 \lesssim \Big(\int_\Omega    \Big( \frac{1}{\delta(z)^2}\int_{B(z, \frac{1}{4}\delta(z))} g^p(w)~dm(w)        \Big)^{q/p} d\mu (z)\Big)^{p/q}
			\end{equation}
		
\noindent By the observation above, the integral in (\ref {ad}) can be bounded by
	
	\begin{equation*}
	\Big(\int_\Omega
\Big(\int_\Omega \frac{1}{\delta(w)^2}
 ~\chi_  {B(w,\frac{1}{3}\delta(w))} (z)~ g^p(w)~dm(w)
       \Big)^{q/p} d\mu (z)\Big)^{p/q}
\end{equation*}		

Since $q/p\geq 1$, by  Minkowski's integral inequality and the condition $(i)$ on the measure, we obtain

\begin{equation*}	
\begin{split}				
\Big(\int_\Omega g^q (z) ~d\mu(z) \Big)^{p/q}\lesssim	
\int_\Omega    \Big(\int_\Omega \Big( \frac{1}{\delta(w)^2}~ \chi_  {B(w,\frac{1}{3}\delta(w))} (z)~ g^p(w)\Big)^{q/p} d\mu(z) \Big)^{p/q} dm(w)\\*[3pt]
\lesssim	\int_\Omega  \frac{1}{\delta(w)^2} ~ g^p(w)~ \big(\mu (B(w,1/3~\delta(w))\big)^{p/q} ~dm(w)\lesssim \int_\Omega  g^p(w)~ \delta(w)^\alpha ~dm(w)
\end{split}		
\end{equation*}	

\noindent with comparison constants only depending on $(C,p,||\mu||)$ as we wanted to prove.

$(ii)\Rightarrow (iii)$
This is an immediate  consequence of the fact that if $f$ is analytic, then $|f|$ is subharmonic.

$(iii)\Rightarrow (i)$
Let $\varphi$ represent a conformal mapping from $\mathbb{D}$ onto $\Omega$. By Koebe's theorem and a change of variables, condition $(iii)$ can be written as
\begin{equation}\label{pull}
\bigg(\int_\mathbb{D}\ |f\circ \varphi|^q  ~d\varphi^\ast (\mu)\bigg)^{1/q}\leq C~ \bigg(\int_\mathbb{D}|f\circ \varphi|^p ~|\varphi'|^{(2+\alpha)} (1-|z|)^\alpha ~dm\bigg)^{1/p}
\end{equation}

\noindent Defining the measure $\tau$ in $\mathbb{D}$ as $d\tau=\frac{1}{|\varphi'|^{(2+\alpha)q/p}}~d\varphi^\ast (\mu)$, (\ref{pull}) is equivalent to
$$
\bigg(\int_\mathbb{D}\ \big(|f\circ \varphi| ~|\varphi'|^{(2+\alpha)/p}\big)^q ~ d\tau \bigg)^{1/q}\leq C~ \bigg(\int_\mathbb{D}\big(|f\circ \varphi| ~|\varphi'|^{(2+\alpha)/p}\big)^p (1-|z|)^\alpha ~dm\bigg)^{1/p}
$$

\noindent which implies by Theorem~\ref{theoA}, that $\tau$ is a $q$-Carleson measure in $\mathbb{D}$ for $\mathcal{A}^p_\alpha$.

Consider now a Whitney ball $B(w_0, R)\subset \Omega$. To show that $\mu(B(w_0,R)\lesssim R^{(2+\alpha)q/p}$ we follow exactly the same steps as in "only if part in Theorem 1, just replacing the exponent $q/p$ by $(2+\alpha) q/p$ and using the fact that since $\tau$ is $q$-Carleson measure in $\mathbb{D}$ for $\mathcal{A}^p_\alpha$, by theorem A, $ \tau(B(z,R)\lesssim R^{(2+\alpha)q/p}$ for all Whitney balls $B(z,R)\subset \mathbb{D}$.

\end{proof}

\textit{Mar\'ia J. Gonz\'alez:} Departamento de Matem\'aticas, Universidad de C\'adiz, Spain. E-mail address: majose.gonzalez@uca.es


\begin{thebibliography}{C-J-S}


\bibitem[A]{A}
 L.V. Ahlfors, Lectures on quasiconformal mappings, Van Nostrand, Princeton
1966; Reprinted by Wadsworth Inc. Belmont, 1987.



\bibitem[BJ]{BJ}
C. J. Bishop; P. W. Jones,
\textit{Harmonic measure, $L^2$ estimates and the Schwarzian derivative.}
J. Anal. Math. 62 (1994), 77-113.


\bibitem[C]{C}
L. Carleson, \textit{Interpolation by bounded analytic functions and the corona problem}, Ann. of Math. (2) 76
(1962), 547-559.

\bibitem[Da]{Da}
G. David, \textit{Courbes corde-arc et espaces de Hardy g\'en\'eralis\'es. (French) [Chord-arc curves and generalized Hardy spaces].}
Ann. Inst. Fourier (Grenoble) 32 (1982), no. 3, xi, 227-239.

\bibitem[D]{D}
P. L. Duren, \textit{Theory of Hp Spaces} Academic Press (New York-London), 1970.(Reprint:
Dover, Mineola, New York 2000).

\bibitem[Du]{Du}
P. L. Duren,
\textit{Extension of a theorem of Carleson.} Bull. Amer. Math. Soc. 75 (1969), 143-146.

\bibitem[G]{G}
J. B. Garnett,  \textit{Bounded analytic functions.}
Academic Press (New York), 1981.


\bibitem[H]{H}
W. W. Hastings, \textit{A Carleson measure theorem for Bergman spaces}, Proc. Amer. Math.
Soc., 52 (1975), 237-241.

\bibitem[HKZ]{HKZ}
H. Hedenmalm, B. Korenblum and K. Zhu, Theory of Bergman Spaces, Graduate
Texts in Mathematics 199, Springer, New York, Berlin, etc. (2000).

\bibitem[JK]{JK}
D. S. Jerison, ; C. E. Kenig,
\textit{Hardy spaces, $A_{\infty }$, and singular integrals on chord-arc domains.}
Math. Scand. 50 (1982), no. 2, 221-247.

\bibitem[L]{L}
D. H. Luecking, \textit{A technique for characterizing Carleson measures on Bergman
spaces}, Proc. Amer. Math. Soc., 87 (1983), 656-660.

\bibitem[Lu]{Lu}
D. H. Luecking,\textit{ Forward and reverse Carleson inequalities for functions in Bergman spaces
and their derivatives}, Amer. J. Math. 107 (1985), 85-111.

\bibitem[OP]{OP}
V. L. Oleinik; B. S. Pavlov, Embedding theorems for weighted classes of harmonic
and analytic functions, (in russian), Zap. Nauch. Sem. LOMI Steklov 22 (1971), 94-
102; english translation in J. Soviet Math. 2 (1974), 135-142.

\bibitem[P]{P}
M. Pavlovi´c, \textit{On subharmonic behaviour and oscillation of functions on balls in $R^n$}, Publ.
Inst. Math. (N.S.) 55(69) (1994), 18–22.


\bibitem[PR]{PR}
M. Pavlovi´c; J. Riihentaus, \textit{Classes of quasi-nearly subharmonic functions}, Potential Analysis
29:1(2008), 89–104

\bibitem[Po]{Po}
Ch. Pommerenke, Boundary behaviour of conformal maps, Springer-Verlag, Berlin, 1992.

\bibitem[S]{S}
D. Stegenga, \textit{Multipliers of the Dirichlet space,} Ill. J. Math., 24 (1980), 113-139.



\bibitem[Z]{Z}
M. Zinsmeister,
\textit{Les domaines de Carleson. (French) [Carleson domains]}
Michigan Math. J. 36 (1989), no. 2, 213-220.

\end{thebibliography}
\end{document}